\documentclass[11pt,a4paper]{amsart}
\usepackage[T1]{fontenc}
\usepackage{lmodern}
\usepackage{amsmath,amssymb,amsthm,mathtools}
\usepackage[a4paper,margin=26mm]{geometry}
\usepackage{microtype,needspace,xcolor}
\definecolor{linkblue}{RGB}{20,55,105}
\usepackage[colorlinks=true,linkcolor=linkblue,citecolor=linkblue,urlcolor=linkblue]{hyperref}
\hypersetup{pdftitle={Small spherical 7-designs in prime dimensions},
 pdfauthor={Andriy Bondarenko and Dmytro Mitin}}
\newtheorem{theorem}{Theorem}
\newtheorem*{corollary}{Corollary}
\newtheorem{proposition}[theorem]{Proposition}
\newtheorem{lemma}[theorem]{Lemma}
\numberwithin{equation}{section}
\title{Small spherical $7$-designs in prime dimensions}
\author{Andriy Bondarenko}
\address{Department of Mathematical Sciences, Norwegian University of
Science and Technology (NTNU), 7491 Trondheim, Norway}
\email{andriybond@gmail.com}
\author{Dmytro Mitin}
\email{dmitin3@gmail.com}
\date{}
\begin{document}
\begin{abstract}
For every prime $p\ge29$, we construct an antipodal
spherical $7$-design on $S^{p-1}$ with $42(p-1)(p+2)^2$ points.
The construction uses eigenlines of finite Weil operators.
\end{abstract}
\maketitle

\section{Introduction}

Let $\sigma$ be normalized surface measure on $S^{d-1}\subset\mathbb R^d$.
A finite set $X\subset S^{d-1}$ is a spherical $t$-design if
\[
 \frac1{|X|}\sum_{x\in X} P(x)=\int_{S^{d-1}}P(x)\,d\sigma(x)
\]
for every polynomial $P$ of degree at most $t$.
Bondarenko, Radchenko and Viazovska \cite{BRV} proved the conjecture
of Korevaar and Meyers that, for fixed $d$,
the minimum number of points in a spherical $t$-design on $S^{d-1}$
is of order $t^{d-1}$ as $t\to\infty$.
Here we fix $t=7$ and let $d$ grow.
The bound of Delsarte, Goethals and Seidel \cite{DGS} gives
$|X|\ge 2\binom{d+2}{3}$ for a spherical $7$-design
\cite[Theorem~1.3]{Dillon}. We prove a bound of the same order in prime dimensions.

\begin{theorem}\label{thm:main}
For every prime $p\ge29$, there is an antipodal
spherical $7$-design $X_p\subset S^{p-1}$ with
\[
 |X_p|=42(p-1)(p+2)^2.
\]
\end{theorem}

To our knowledge, even for positive
quadrature it was open whether a formula of degree 7 with $O(d^3)$
nodes exists along an unbounded sequence of dimensions.
Kuperberg \cite[Corollary~4.2 and the remark following it]{Kuperberg}
obtained positive formulas of degree seven with $O(d^4)$ nodes in
every dimension and equal weights when $d=4\cdot7^k$.
Dillon \cite[Theorem~D]{Dillon} obtained $O(d^3)$ nodes in every
dimension with weights of either sign. The following corollary answers
the case $t=3$ of his question on positive weights \cite[Section~7, Question~1]{Dillon}.

\begin{corollary}
For every $d\ge1$, there is a positive quadrature formula of degree $7$
on $S^{d-1}$ with $O(d^3)$ nodes.
\end{corollary}

\begin{proof}
Applying Dillon's projection argument \cite[Corollary~3.5]{Dillon} to
$X_p$ for the least prime $p\ge\max\{29,d\}$ gives a positive formula
on $S^{d-1}$ of degree $7$ with at most $16|X_p|=O(d^3)$ nodes.
\end{proof}

For a unit vector $y$, the first six moments of $\langle x,y\rangle$
on $S^{d-1}$, in increasing order, are
\[
 0,\quad \frac1d,\quad 0,\quad \frac3{d(d+2)},\quad 0,\quad
 \frac{15}{d(d+2)(d+4)}.
\]
An antipodal set $X$ is a $7$-design if and only if it reproduces
the sixth moment in every direction, that is,
\begin{equation}\label{eq:sixth-criterion}
 \frac1{|X|}\sum_{x\in X}\langle x,y\rangle^6
 =\frac{15}{d(d+2)(d+4)}
 \qquad (y\in S^{d-1}).
\end{equation}
By homogeneity this determines the sixth moment for every $y\in\mathbb R^d$.
Comparing coefficients in $y$ then gives the required equality for every
monomial of degree 6. For a homogeneous polynomial of degree 2 or 4,
multiplication by a power of $\|x\|^2$ raises its degree to 6 and
preserves its values on the sphere. All odd moments through degree 7 vanish by antipodality.

Recently, Arman, Bondarenko, Prymak and Radchenko \cite{ABPR}
lifted complex projective designs to spherical 5-designs with $O(d^2)$
points in every dimension. For each maximal torus of $\mathrm{SL}_2(\mathbb F_p)$,
Gurevich, Hadani and Sochen \cite[Section~IV]{GHS} chose a unit vector
in each one-dimensional common eigenspace of its Weil operators.
They called the resulting collection the finite oscillator system and
proved bounds of order $p^{-1/2}$ for its cross-correlations under time
and frequency shifts. Following their terminology, we call a torus
\emph{split} if it is diagonalizable over $\mathbb F_p$ and \emph{nonsplit}
if the quadratic extension $\mathbb F_{p^2}$ is needed. We use the odd
eigenfunctions from both families and a third family indexed by
nilpotent matrices.

For unit $z\in\mathbb C^{(p-1)/2}$, all three families of eigenvectors $u$
have the same averages of $|\langle z,u\rangle|^2$ and $|\langle z,u\rangle|^4$. For the sixth powers we need a particular linear combination
of the three sums. We obtain its coefficients by choosing the last
coordinates of our points.

\Needspace{9\baselineskip}
In Section~\ref{sec:construction} we define $X_p$ by placing regular
polygons on the eigenlines and derive the equations for their squared heights from the
moment identities in Proposition~\ref{prop:angular-moments}.
Lemma~\ref{lem:valid-heights} gives the heights needed for
\eqref{eq:sixth-criterion}. Section~\ref{sec:families} describes the
three families of eigenlines of the Weil operators. In
Section~\ref{sec:angular} we prove Proposition~\ref{prop:angular-moments}
from the matrix entries in~\eqref{eq:explicit-U}.

\medskip
\noindent\textbf{Acknowledgments.}
A. Bondarenko received partial support from the Research Council of
Norway through Grant 334466. The authors gratefully acknowledge the
assistance of OpenAI's ChatGPT, whose exploratory discussions, calculations and
numerics helped stimulate the development of this paper.

\Needspace{16\baselineskip}
\section{The points of the design}\label{sec:construction}

Fix a prime $p\ge29$. Write $\mathbb F_p=\mathbb Z/p\mathbb Z$ and
$m=(p-1)/2$. We write each point of
$\mathbb R^p=\mathbb C^m\oplus\mathbb R$ as $x=(w,x_p)$, where
$x_p$ is its last real coordinate. The points of $X_p$ will have
$x_p^2\in\{h_1,\ldots,h_8\}$. These squared heights are determined
by \eqref{eq:constant-heights} and Lemma~\ref{lem:valid-heights}.
In particular, $h_1=h_2=h_3=5/(6p)$.

To specify the complex eigenlines, take the nonzero matrices
\[
 X=\begin{pmatrix}r&s\\t&-r\end{pmatrix}
 \quad (r,s,t\in\mathbb F_p)
\]
up to multiplication by a nonzero scalar. The matrices in
$\mathrm{SL}_2(\mathbb F_p)$ commuting with $X$ form a cyclic group
$C(X)$ which depends only on the scalar class of $X$.
For any generator $g_X$ of this group, the Weil operator $U_{g_X}$ in
\eqref{eq:explicit-U} has $m$ mutually orthogonal eigenlines on the
space of odd functions $f$ on $\mathbb F_p$. The eigenlines are independent
of the generator (Lemma~\ref{lem:groups-bases}). We regard them as lines in
$\mathbb C^m$ by identifying $f$ with $\sqrt2\,(f(1),\ldots,f(m))$.

The type in Table~\ref{tab:families} depends only on the scalar class,
since scaling $X$ multiplies $-\det X$ by a square. The table gives
the number $b_a$ of classes of each type and $c_a=|C(X)|/2$.

\begin{table}[!ht]
\caption{Vertex counts for one eigenline in one eigenbasis}\label{tab:families}
\centering
\setlength{\arraycolsep}{4pt}
$
\begin{array}{c|c|c|c|c|c}
 a&-\det X&b_a&c_a&\text{squared height }h&\text{count for }N(L,h)\\ \hline
 s\ (\text{split}) &\text{nonzero square}&p(p+1)/2&(p-1)/2&h_1=h_2=h_3&21\\
 &&&&h_4,\ h_5,\ h_6&7\text{ each}\\ \hline
 n\ (\text{nonsplit})&\text{nonsquare}&p(p-1)/2&(p+1)/2&h_7&42\\
 0\ (\text{nilpotent})&0&p+1&p&h_8&168
\end{array}$
\end{table}

Fix an eigenline $L$ and a squared height $h\in\{h_1,\ldots,h_8\}$.
We write $N(L,h)$ for the number of vertices at last coordinate $\sqrt h$
whose complex component lies in $L$. Read the \emph{rightmost column
of Table~\ref{tab:families}} once for each matrix class whose eigenbasis
contains $L$, using the row for its type and squared height $h$.
If a single entry applies, $N(L,h)$ is the number printed there. If
several entries apply, $N(L,h)$ is their sum. Put $N(L,h)=0$ if there
are no entries. The $21$ at
$h_1=h_2=h_3$ is counted once per split basis containing $L$.

Let $u_L$ be the unit vector on $L$ whose first nonzero coordinate is
real and positive. At last coordinate $\sqrt h$ we place a regular polygon
with the prescribed $N(L,h)$ vertices whenever $N(L,h)>0$. The union
of these polygons and their antipodes is
\begin{equation}\label{eq:angles}
 X_p=\bigcup_{L,h:\,N(L,h)>0}
 \left\{\,\pm\left(\sqrt{1-h}\,e^{2\pi i k/N(L,h)}u_L,\sqrt h\right)
              \ \middle|\ k=0,\ldots,N(L,h)-1\,\right\}.
\end{equation}
All the vertices are distinct. Indeed, $0<h<1$, so the last coordinate
determines the sign and $h$, and the nonzero complex component determines
$L$ and then $k$. Summing the prescribed vertex counts gives
\begin{equation}\label{eq:point-count}
 K:=\frac17\sum_{L,h}N(L,h)
 =m(6b_s+6b_n+24b_0)=3(p-1)(p+2)^2,\qquad |X_p|=14K.
\end{equation}

\subsection{The equations for the heights}\label{subsec:heights}
For each matrix class $X$, write $u_{X,1},\ldots,u_{X,m}$ for the vectors
$u_L$ on its $m$ eigenlines, in any order. For $a=s,n,0$, let $\mathcal U_a$ consist
of these vectors for all $X$ of type $a$, counted once per basis.
Section~\ref{sec:angular} proves the following identities.

\begin{proposition}\label{prop:angular-moments}
For $z$ on the unit sphere $S^{2m-1}=S^{p-2}$ of $\mathbb C^m$ and each type $a=s,n,0$,
\begin{align}
 \sum_{u\in\mathcal U_a}|\langle z,u\rangle|^2
 &= b_a,\label{eq:angular-second}\\
 \sum_{u\in\mathcal U_a}|\langle z,u\rangle|^4
 &=\frac{2b_a}{m+1}.\label{eq:angular-fourth}
\end{align}
For the sixth powers,
\begin{equation}\label{eq:angular-sixth}
 \sum_{a\in\{s,n,0\}}c_a^2\sum_{u\in\mathcal U_a}|\langle z,u\rangle|^6
 =6p^2.
\end{equation}
\end{proposition}

Write $y\in S^{p-1}$ as $y=(s z,r)$ with $z\in S^{p-2}$, $s\ge0$ and
$s^2+r^2=1$. At a pole, $s=0$ and $z$ may be chosen arbitrarily.
The real inner product of $y$ with
$(\sqrt{1-h}\,e^{i\theta}u,\sqrt h)$ is
$r\sqrt h+s\sqrt{1-h}\,\operatorname{Re}(e^{i\theta}\langle z,u\rangle)$.
The sixth power has the same average on each of our polygons as on
its circle. Indeed, its trigonometric degree is at most $6$, and for
$N=N(L,h)\ge7$ we have $\sum_{k=0}^{N-1} e^{2\pi i k j/N}=0$
for $1\le |j|\le6$.
The averages of $\cos^2\theta,\cos^4\theta,\cos^6\theta$ are
$1/2,3/8,5/16$. Thus the sixth power averages to
\begin{equation}\label{eq:circle-sixth}
\begin{aligned}
 &r^6h^3
 +\frac{15}{2}r^4s^2 h^2(1-h)|\langle z,u\rangle|^2\\
 &\qquad
 +\frac{45}{8}r^2s^4 h(1-h)^2|\langle z,u\rangle|^4
 +\frac5{16}s^6(1-h)^3|\langle z,u\rangle|^6 .
\end{aligned}
\end{equation}
The antipodal polygon has the same sixth moment. Together the two
polygons have weight $2N(L,h)/|X_p|=N(L,h)/(7K)$ in the average over $X_p$.

At every vertex in \eqref{eq:angles} we have $x_p^2=h$. The moment
formulas in the introduction, applied to the last coordinate vector
$e_p=(0,1)\in\mathbb C^m\oplus\mathbb R$, give the spherical averages
of $x_p^2,x_p^4,x_p^6$ as
$1/p$, $3/[p(p+2)]$ and $15/[p(p+2)(p+4)]$, respectively.
In the split family each of $h_1,\ldots,h_6$ contributes weight $mb_s/K$,
with the first three weights added at their common value. The weights
for $h_7$ and $h_8$ are $6mb_n/K$ and $24mb_0/K$, respectively.
Thus the first three moments of $x_p^2$ agree with their spherical
averages precisely when, for every polynomial $f$ of degree at most three,
\begin{equation}\label{eq:height-targets}
\begin{aligned}
 \frac{p(p+1)}{12}\sum_{i=1}^6 f(h_i)
 +\frac{p(p-1)}2 f(h_7)+4(p+1)f(h_8)
   =(p+2)^2\!\int_{S^{p-1}} f(x_p^2)\,d\sigma(x).
\end{aligned}
\end{equation}

At the equator $r=0$ and $s=1$ in \eqref{eq:circle-sixth}. The coefficients of the
three sums of $|\langle z,u\rangle|^6$, apart from the common
factor $5/(16K)$, are
\[
 \sum_{i=1}^6(1-h_i)^3,\qquad 6(1-h_7)^3,\qquad 24(1-h_8)^3.
\]
To apply \eqref{eq:angular-sixth}, we choose them in the ratio
$c_s^2:c_n^2:c_0^2$. After division by $6,6,24$, respectively,
the ratios become $(1-1/p)^2:(1+1/p)^2:1$. The average of these
last three numbers with the family weights $6mb_s/K$, $6mb_n/K$ and $24mb_0/K$ is
$(p+1)(p+3)/(p+2)^2$ by Table~\ref{tab:families}. On the other hand,
\eqref{eq:height-targets} prescribes the average of $(1-h_i)^3$ as
\[
 \int_{S^{p-1}}(1-x_p^2)^3\,d\sigma(x)
 =\frac{(p-1)(p+1)(p+3)}{p(p+2)(p+4)}.
\]
Dividing the two averages gives the factor in
\begin{equation}\label{eq:height-sixth}
 \left(\frac16\sum_{i=1}^6(1-h_i)^3,\ (1-h_7)^3,\ (1-h_8)^3\right)
 =\frac{(p-1)(p+2)}{p(p+4)}
 \left(\left(1-\frac1p\right)^2,\left(1+\frac1p\right)^2,1\right).
\end{equation}
The last two entries determine $h_7,h_8$. The first entry then
follows from \eqref{eq:height-targets} with $f(t)=(1-t)^3$.

\subsection{Choice of the heights}
Take
\begin{equation}\label{eq:constant-heights}
\begin{aligned}
 h_1=h_2=h_3&=\frac5{6p},\\
 (h_7,h_8)&=(1,1)-
 \sqrt[3]{\frac{(p-1)(p+2)}{p(p+4)}}\,
 \bigl((1+1/p)^{2/3},1\bigr).
\end{aligned}
\end{equation}
The three moment equations \eqref{eq:height-targets} determine the
remaining squared heights $h_4,h_5,h_6$.

\begin{lemma}\label{lem:valid-heights}
For $p\ge29$, equations \eqref{eq:height-targets} and \eqref{eq:constant-heights}
determine unique $h_4,h_5,h_6$ with
\[
 0<ph_4<1<ph_5<3<ph_6<6.
\]
All eight squared heights lie in $(0,1)$ and satisfy
\eqref{eq:height-sixth}.
\end{lemma}

\begin{proof}
Write $a=ph_4$, $b=ph_5$ and $c=ph_6$. We seek three numbers with
prescribed sum, sum of squares and sum of cubes,
\[
 a+b+c=s_1,\qquad a^2+b^2+c^2=s_2,\qquad a^3+b^3+c^3=s_3.
\]
To compute the right-hand sides, apply \eqref{eq:height-targets} to
$f(t)=pt$, $f(t)=(pt)^2$ and $f(t)=(pt)^3$, divide by $p(p+1)/12$ and
subtract the contributions of the five heights fixed in \eqref{eq:constant-heights}.
This gives
\[
\begin{aligned}
 s_1&=\frac{12(p+2)^2}{p(p+1)}-\frac52
       -\frac{6(p-1)}{p+1}(ph_7)-\frac{48}{p}(ph_8),\\
 s_2&=\frac{36(p+2)}{p+1}-\frac{25}{12}
       -\frac{6(p-1)}{p+1}(ph_7)^2-\frac{48}{p}(ph_8)^2,\\
 s_3&=\frac{180p(p+2)}{(p+1)(p+4)}-\frac{125}{72}
       -\frac{6(p-1)}{p+1}(ph_7)^3-\frac{48}{p}(ph_8)^3.
\end{aligned}
\]
Expanding $(a+b+c)^2$ and $(a+b+c)^3$ gives
\[
 ab+bc+ca=\frac{s_1^2-s_2}{2},\qquad
 abc=\frac{s_1^3-3s_1s_2+2s_3}{6}.
\]
Thus $(t-a)(t-b)(t-c)$ must be the polynomial
\[
 F_p(t)=t^3-s_1t^2+\frac{s_1^2-s_2}{2}\,t
                 -\frac{s_1^3-3s_1s_2+2s_3}{6}.
\]
Conversely, the same formulas show that the roots of $F_p$ have power
sums $s_1,s_2,s_3$. It remains to locate the roots.

By \eqref{eq:constant-heights}, $ph_7\to 1/3$ and $ph_8\to1$ as
$p\to\infty$. The formulas above therefore give
\[
 (s_1,s_2,s_3)\longrightarrow
 \left(\frac{15}{2},\frac{133}{4},\frac{4273}{24}\right).
\]
The choice $h_1=h_2=h_3=5/(6p)$ now gives a limiting polynomial
$36t^3-270t^2+414t-179$ for $36F_p(t)$. Its values at $0,1,3,6$
are $-179,1,-395,361$. For $p\ge29$, $F_p$ has the same signs at these four points.
For $p>50$ this follows by estimating the difference between $36F_p$ and
its limiting cubic at these points. For $p=29,31,37,41,43,47$, the signs
can be checked by direct numerical substitution.
Hence $F_p$ has one root in each of $(0,1)$, $(1,3)$, and $(3,6)$.
Take them in order as $ph_4,ph_5,ph_6$. This gives the unique heights
in the stated intervals. Formula~\eqref{eq:constant-heights} also gives
$0<h_7<h_8<1$ by cubing. Finally, \eqref{eq:height-targets} with
$f(t)=(1-t)^3$ gives the remaining entry of \eqref{eq:height-sixth}.
\end{proof}

\begin{proof}[Proof of Theorem~\ref{thm:main}]
Use these heights in \eqref{eq:angles}. The set is antipodal, and
\eqref{eq:point-count} gives $42(p-1)(p+2)^2$ distinct points.
It remains to average \eqref{eq:circle-sixth} with weights $N(L,h)/(7K)$.
By \eqref{eq:height-targets}, the averages of $h^3$, $h^2(1-h)$,
and $h(1-h)^2$ are, respectively,
\[
 \frac{15}{p(p+2)(p+4)},\qquad
 \frac{3(p-1)}{p(p+2)(p+4)},\qquad
 \frac{p^2-1}{p(p+2)(p+4)}.
\]
Each family $\mathcal U_a$ has $mb_a$ members, so
\eqref{eq:angular-second} and \eqref{eq:angular-fourth} give family
averages $1/m$ and $2/[m(m+1)]$. Using $2m=p-1$, the first three terms
of \eqref{eq:circle-sixth} therefore contribute
\[
 \frac{15}{p(p+2)(p+4)}\bigl(r^6+3r^4s^2+3r^2s^4\bigr).
\]
For the last term, \eqref{eq:height-sixth} and \eqref{eq:angular-sixth}
give the coefficient of $s^6$ as
\[
 \frac5{16K}\,\frac{24(p-1)(p+2)}{p^3(p+4)}\,6p^2
 =\frac{15}{p(p+2)(p+4)}.
\]
The four terms thus give $15(r^2+s^2)^3/[p(p+2)(p+4)]$.
Since $r^2+s^2=1$, this proves \eqref{eq:sixth-criterion}.
\end{proof}

\section{The Weil operators and their eigenbases}\label{sec:families}

\subsection{The operators}\label{subsec:explicit-vectors}

Identify $\mathbb C^m$ with the space $H$ of odd functions
$f:\mathbb F_p\to\mathbb C$, with inner product
$\langle f,g\rangle=\sum_x\overline{f(x)}g(x)$ and orthonormal basis
\[
 \frac{\delta_j-\delta_{-j}}{\sqrt2},\qquad 1\le j\le m.
\]
Thus $f\in H$ has coordinates $(\sqrt2\,f(1),\ldots,\sqrt2\,f(m))$.
We use $M^*$ for the conjugate transpose of a complex matrix $M$.

On all functions on $\mathbb F_p$, define
\[
 (D_{r,s}f)(x)=
 \exp\!\left(\frac{2\pi i}{p}\left(sx+\frac{rs}{2}\right)\right)f(x+r).
\]
In each exponential, the coefficient of $2\pi i/p$ is computed
in $\mathbb F_p$ and represented by any integer.
For $g=\left(\begin{smallmatrix}a&b\\c&d\end{smallmatrix}\right)
\in\mathrm{SL}_2(\mathbb F_p)$,
the Weil operator $U_g$ acts on functions and carries the index
$v=(r,s)$ of $D_v$ to the column vector $gv=(ar+bs,cr+ds)$. Precisely,
\begin{equation}\label{eq:intertwining}
 U_gD_v=D_{gv}U_g\qquad(v\in\mathbb F_p^2).
\end{equation}
Thus applying $D_v$ before $U_g$ has the same effect as applying
$D_{gv}$ afterwards. This identity determines the unitary $U_g$
up to a scalar of absolute value one \cite[Section II-E]{GHS}.
We use the choice
\begin{equation}\label{eq:explicit-U}
 (U_gf)(x)=
 \begin{cases}
 \displaystyle\frac1{\sqrt p}\sum_{y\in\mathbb F_p}
 \exp\!\left(\frac{2\pi i}{p}
              \frac{-dx^2+2xy-ay^2}{2b}\right)f(y),&b\ne0,\\[7pt]
 \displaystyle\exp\!\left(-\frac{2\pi i}{p}
              \frac{cx^2}{2a}\right)f(a^{-1}x),&b=0.
 \end{cases}
\end{equation}
For example, when $g=\bigl(\begin{smallmatrix}0&1\\-1&0\end{smallmatrix}\bigr)$,
$U_g$ is the finite Fourier transform with exponential $e^{2\pi ixy/p}$.
Replacing $f(y)$ by $f(y+r)$ in its sum and setting $y'=y+r$ gives
the factor $e^{-2\pi irx/p}$. In other words,
$U_gD_{r,0}=D_{0,-r}U_g$, in agreement with $g(r,0)=(0,-r)$.
For general $g$, substitution in \eqref{eq:explicit-U} gives
$U_gD_{r,0}=D_{ar,cr}U_g$ and $U_gD_{0,s}=D_{bs,ds}U_g$.
Multiplying the last two identities gives \eqref{eq:intertwining}:
since $ad-bc=1$, the scalar factor in $D_{ar,cr}D_{bs,ds}$ is the
same as in $D_{r,0}D_{0,s}$. For $b\ne0$, \eqref{eq:explicit-U} is a composition of
a finite Fourier transform, a change of variable and multiplications
by quadratic exponentials. Each operation is unitary and preserves
odd functions. The same is true of the formula for $b=0$.

Two successive changes of index give
\[
 (U_gU_h)D_v=U_gD_{hv}U_h=D_{ghv}(U_gU_h).
\]
Since $U_{gh}$ also carries $D_v$ to $D_{ghv}$, the operator
$U_{gh}^*U_gU_h$ commutes with every $D_v$. The $D_v$ span all
complex $p\times p$ matrices. Indeed, for fixed $r$ their entries in positions
$(x,x+r)$ form a Fourier basis as $s$ varies, and varying $r$ covers
every matrix entry. It follows that $U_{gh}^*U_gU_h$ is scalar.
Thus $U_gU_h$ differs from $U_{gh}$ only by a scalar factor, which
cancels in conjugation.

\subsection{The eigenbases}

Represent the $p^2+p+1$ classes of nonzero matrices with trace $0$ by
\[
 X=\begin{pmatrix}r&s\\t&-r\end{pmatrix},\qquad
 (r,s,t)=(1,s,t),\ (0,1,t),\ \text{or }(0,0,1),
\]
where the free entries range over $\mathbb F_p$.
Put $\Delta=-\det X$, so that $X^2=\Delta I$. The centralizer of $X$ in
$\mathrm{SL}_2(\mathbb F_p)$ is
\begin{equation}\label{eq:explicit-centralizer}
 C(X)=\{aI+bX:a,b\in\mathbb F_p,\ a^2-\Delta b^2=1\}.
\end{equation}

\begin{lemma}[Groups and eigenbases]\label{lem:groups-bases}
Each group $C(X)$ is cyclic. Its order and the number of classes
of each type are given in Table~\ref{tab:families}.
If $g_X$ generates $C(X)$, then $U_{g_X}$ has $m$ distinct eigenvalues
on $H$. Any two generators of $C(X)$ give the same $m$ eigenlines.
\end{lemma}

\begin{proof}
Since $X$ is not a scalar matrix, some
$v\in\mathbb F_p^2$ has $v,Xv$ linearly independent.  A matrix
commuting with $X$ is determined by its value on $v$.  If this
value is $av+bXv$, the matrix is $aI+bX$.  Its determinant is
$a^2-\Delta b^2$, proving the formula for $C(X)$.

If $\Delta$ is a nonzero square, choose an eigenbasis of $X$
whose change-of-basis matrix has determinant one.
The elements of $C(X)$ then become
$\operatorname{diag}(a,a^{-1})$, $a\in\mathbb F_p^\times$.
Thus $C(X)$ is cyclic of order $p-1$.
For a generator $a$ of $\mathbb F_p^\times$, the corresponding
operator on $H$ is, up to a scalar and a unitary change of basis,
$f(x)\mapsto f(a^{-1}x)$. It permutes the $m$ lines
$\mathbb C(\delta_j-\delta_{-j})$ in one cycle, and its $m$th power is $-I$.
The scaling operator has characteristic polynomial $t^m+1$,
so its eigenvalues, and hence those of $U_{g_X}$, are distinct.

If $\Delta$ is a nonsquare, each nonzero matrix $aI+bX$ has nonzero
determinant $a^2-\Delta b^2$. These matrices, together with zero,
form a field with $p^2$ elements in which
\[
 (aI+bX)^p=aI-bX,
 \qquad (aI+bX)^{p+1}=(a^2-\Delta b^2)I.
\]
Among the $p^2-1$ elements of its cyclic multiplicative group, the
$p+1$ solutions of $g^{p+1}=I$ form $C(X)$.
To prove distinctness of the eigenvalues, first consider $U_{g_X}$
on all functions on $\mathbb F_p$.  Conjugation by $U_{g_X}$
permutes the basis $D_v$ by $v\mapsto g_Xv$.
Every nonzero orbit has $p+1$ elements, since $g_X^k-I$ is an
invertible matrix whenever $g_X^k\ne I$.  There are therefore
$1+(p^2-1)/(p+1)=p$ orbits.  A matrix commutes with $U_{g_X}$
exactly when its coefficients in the basis $D_v$ are constant on
these orbits, so the space of commuting matrices has dimension $p$.
If the eigenspaces of $U_{g_X}$ have dimensions $d_1,d_2,\ldots$,
that same dimension is $\sum_jd_j^2$.  Since $\sum_jd_j=p$,
every $d_j$ is one.  This also proves distinctness on $H$.

If $\Delta=0$, then
$C(X)=\{\pm(I+bX):b\in\mathbb F_p\}$.
The element $-(I+X)$ has order $2p$ and generates this group.
Every generator is $-(I+bX)$ with $b\ne0$.
Choose a basis of $\mathbb F_p^2$ of determinant one whose second
vector spans $\ker X$. Then $g_X$ is lower triangular with both
diagonal entries $-1$. On odd functions, \eqref{eq:explicit-U}
gives, up to a scalar and unitary conjugation, multiplication by
\[
 \exp\!\left(\frac{2\pi i}{p}\frac{tx^2}{2}\right),
 \qquad t\ne0.
\]
Its eigenvalues on $\delta_j-\delta_{-j}$, $1\le j\le m$, are distinct,
because $x^2=y^2$ implies $x=\pm y$.

If $g$ and $g'$ both generate $C(X)$, each is a power of the
other. By the product formula above, $U_{g'}$ is a scalar multiple
of a power of $U_g$, so every eigenvector of $U_g$ is an eigenvector
of $U_{g'}$. Interchanging $g$ and $g'$ gives the reverse inclusion.
Since the eigenspaces on $H$ are one-dimensional, the two operators
give the same lines.

Finally, the split classes correspond to unordered pairs of distinct
eigenlines in $\mathbb F_p^2$, of which there are
$\binom{p+1}{2}=p(p+1)/2$. A nonzero matrix with $X^2=0$ is determined
up to scale by its kernel, giving $p+1$ classes. The remaining
$p(p-1)/2$ classes are nonsplit.
\end{proof}

\section{Moments of the complex vectors}\label{sec:angular}

\subsection{Averages in an eigenbasis}
For each class $X$, put $c_X=|C(X)|/2$. Each $g\in C(X)$ is a
$2\times2$ matrix of determinant one over $\mathbb F_p$ and acts on
$\mathbb F_p^2$. Its Weil operator $U_g$ preserves the complex space
$H$ of odd functions. In the basis of Section~\ref{subsec:explicit-vectors},
the restriction of $U_g$ to $H$ is a unitary $m\times m$ matrix,
where $m=(p-1)/2$, so $U_g^*=U_g^{-1}$. A complex-linear map
$A:H\to H$ is represented in this basis by an arbitrary complex
$m\times m$ matrix. We average $A$ over $C(X)$ by setting
\[
 E_X(A)=\frac{1}{2c_X}\sum_{g\in C(X)}U_g A U_g^*.
\]
In the eigenbasis $u_{X,1},\ldots,u_{X,m}$, this average keeps the
diagonal of $A$ and sets its other entries to zero. Indeed, every
$U_g$ in the sum preserves the eigenlines, so conjugation leaves the
diagonal entries unchanged. Conjugating the whole sum by $U_{g_X}$
permutes its terms. Hence $E_X(A)$ commutes with $U_{g_X}$, whose
eigenvalues are distinct by Lemma~\ref{lem:groups-bases}.

\begin{lemma}\label{lem:matrix-averages}
For any two complex $m\times m$ matrices $A,B$ and each type $a=s,n,0$,
\begin{equation}\label{eq:second-group-sum}
 \sum_{X\text{ of type }a}E_X(A)
 =\frac{b_a}{m+1}\bigl(A+\operatorname{tr}(A)I\bigr).
\end{equation}
The sum over all three types satisfies
\begin{equation}\label{eq:third-group-sum}
 \sum_X c_X^2 E_X(A)E_X(B)
 =p^2\bigl[(A+\operatorname{tr}(A)I)(B+\operatorname{tr}(B)I)
          +BA+\operatorname{tr}(AB)I\bigr].
\end{equation}
\end{lemma}

\begin{proof}[Proof of Proposition~\ref{prop:angular-moments}]
For a fixed class $X$, expand $z$ in the orthonormal eigenbasis of the Weil
operator $U_{g_X}$,
\[
 z=\sum_{j=1}^m \alpha_j u_{X,j},\qquad
 \alpha_j=\langle u_{X,j},z\rangle .
\]
Since $z$ is a unit vector, $\sum_j|\alpha_j|^2=1$. The sum over $\mathcal U_a$ counts
each eigenvector once for every basis in which it occurs. Summing over
the $b_a$ bases of type $a$ therefore gives \eqref{eq:angular-second}.

For the fourth and sixth powers, let $A=zz^*$,
where $z$ is a column vector and $z^*=(\overline{z_1},\ldots,\overline{z_m})$
is its conjugate transpose. Since $Aw=\langle z,w\rangle z$, $A$ is the
orthogonal projection onto $\mathbb Cz$. In the eigenbasis for $X$, its
$(i,j)$ entry is $\alpha_i\overline{\alpha_j}$. Thus $E_X(A)$ has
diagonal entries $|\alpha_j|^2$ and all its other entries are zero.
Squaring $E_X(A)$ squares each diagonal entry. Taking the inner product
with $z$ supplies one more factor $|\alpha_j|^2$ in each summand, giving
\[
 \begin{aligned}
 z^*E_X(A)z
   &=\sum_{j=1}^m\overline{\alpha_j}\,|\alpha_j|^2\alpha_j
     =\sum_{j=1}^m|\alpha_j|^4,\\
 z^*E_X(A)^2z
   &=\sum_{j=1}^m\overline{\alpha_j}\,|\alpha_j|^4\alpha_j
     =\sum_{j=1}^m|\alpha_j|^6 .
 \end{aligned}
\]

Summing the first line over the classes $X$ of type $a$ and using
\eqref{eq:second-group-sum} and $\operatorname{tr}(A)=1$, we obtain
\[
 \sum_{u\in\mathcal U_a}|\langle z,u\rangle|^4
 =z^*\Bigl(\sum_{X\text{ of type }a}E_X(A)\Bigr)z
 =\frac{b_a}{m+1}\,z^*(A+I)z
 =\frac{2b_a}{m+1},
\]
since $Az=z$.

For the sixth powers, multiply the second line by $c_X^2$ and sum
over all $X$. Grouping the terms by type gives the weights $c_a^2$ in
\eqref{eq:angular-sixth}, since $c_X=c_a$ for $X$ of type $a$.
Now take $B=A$ in \eqref{eq:third-group-sum}. The three terms on the
right involve $(A+I)^2$, $A^2$ and $\operatorname{tr}(A^2)I$.
We have $(A+I)z=2z$, $A^2=A$ and $\operatorname{tr}(A^2)=1$. Hence
\[
 \sum_{a\in\{s,n,0\}}c_a^2\sum_{u\in\mathcal U_a}|\langle z,u\rangle|^6
 =z^*\Bigl(\sum_X c_X^2 E_X(A)^2\Bigr)z
 =p^2(4+1+1)=6p^2 .\qedhere
\]
\end{proof}

\subsection{Sums of Weil matrices}

\begin{proof}[Proof of Lemma~\ref{lem:matrix-averages}]
Consider the determinant-one matrices with equal diagonal entries,
\[
 \mathcal M=\left\{
 M=\begin{pmatrix}a&b\\c&a\end{pmatrix}:
 a,b,c\in\mathbb F_p,\quad a^2-bc=1\right\}.
\]
We compare the sum over all $M$ with the sum over those for which $XM$
also has equal diagonal entries. If $b\ne0$, the entries $a,b$ determine $c$.
If $b=0$, we have $a=\pm1$ and $c$ is arbitrary. Thus $\mathcal M$ has
$p(p-1)+2p=p(p+1)$ members, with $M$ and $-M$ counted separately.
Let $A^{\mathsf T}$ denote the ordinary transpose in the basis
$(\delta_j-\delta_{-j})/\sqrt2$ of $H$. For each $M$ form the $m\times m$ matrix
\[
 T_M(A)=U_M A^{\mathsf T}U_M^* .
\]

\smallskip
\noindent\emph{The sum over all $M$.}
We first show that
\begin{equation}\label{eq:matrix-total}
 \sum_{M\in\mathcal M} T_M(A)=2p\bigl(A+\operatorname{tr}(A)I\bigr).
\end{equation}
Write $\omega=e^{2\pi i/p}$. When $b\ne0$ put $h=-a/(2b)$ and $k=1/b$.
Then $h$ runs through $\mathbb F_p$ and $k$ through $\mathbb F_p^\times$ independently.
In our basis of odd functions, \eqref{eq:explicit-U} gives
\[
 (U_M)_{xy}=\frac{1}{\sqrt p}\,\omega^{h(x^2+y^2)}
                 (\omega^{kxy}-\omega^{-kxy}),\qquad 1\le x,y\le m.
\]
Thus $U_M$ is symmetric. For $b=0$ it is diagonal. For each $a=\pm1$,
as $c$ varies, we obtain all $p$ matrices $\operatorname{diag}(\omega^{h x^2})$
up to an overall sign, with $h=-c/(2a)$.

By linearity we may take $A$ to have a single $1$ in row $i$, column $j$,
and zeros elsewhere. Then $(T_M(A))_{xy}=(U_M)_{xj}\overline{(U_M)_{yi}}$.
For $b\ne0$, sum this product first over $h$ and then over $k$. We may
include $k=0$, since the displayed difference vanishes there.
Using $\sum_{h\in\mathbb F_p}\omega^{hd}=0$ for $d\ne0$, we get zero unless
\[
 x^2+j^2=y^2+i^2,\qquad xj=\pm yi .
\]
Thus $x^2,j^2$ have the same sum and product as $y^2,i^2$.
Since $1,\ldots,m$ have distinct squares, either $x=y$ and $j=i$,
or $x=i$ and $j=y$. In either case $(U_M)_{xj}=(U_M)_{yi}$.
The sum over $h$ cancels the factor $1/p$, leaving
\[
 \sum_{k\in\mathbb F_p}|\omega^{kxj}-\omega^{-kxj}|^2
  =\sum_{k\in\mathbb F_p}(2-\omega^{2kxj}-\omega^{-2kxj})=2p ,
\]
since $2xj\ne0$. For $b=0$ the product vanishes unless $x=j$ and $y=i$.
In that case it is $\omega^{h(x^2-y^2)}$. Summing over $h$ and both choices
of $a$ gives a further $2p$ precisely when $x=y=i=j$, and zero otherwise.
Hence for $i\ne j$ the sum is $2pA$. For $i=j$ it is $2p(I+A)$.
This proves \eqref{eq:matrix-total}.

\smallskip
\noindent\emph{The sum for a fixed $X$.}
Write $X=\left(\begin{smallmatrix}r&s\\t&-r\end{smallmatrix}\right)$ and
$\Delta=-\det X=r^2+st$. The diagonal entries of $XM$ are $ra+sc$ and $tb-ra$.
We claim that the matrices for which these entries agree give
\begin{equation}\label{eq:matrix-local}
  \sum_{\substack{M\in\mathcal M\\2ar+cs-bt=0}}T_M(A)
     = 2c_X E_X(A).
\end{equation}
Let $U=U_{g_X}$, where $g_X$ generates $C(X)$. Recall that $E_X(A)$
keeps the diagonal of $A$ in the eigenbasis of $U$. If $M$ occurs in the
sum and $g\in C(X)$, then $gM$ also occurs. Indeed, $gM$ and $XgM$
are linear combinations of $M$ and $XM$, since $g$ has the form
$uI+vX$ and $X^2=\Delta I$. Thus $M\mapsto gM$ permutes the terms. As
$T_{gM}(A)=U_g T_M(A)U_g^*$, the sum commutes with $U$.
By Lemma~\ref{lem:groups-bases}, $U$ has distinct eigenvalues, so the
sum is diagonal in its eigenbasis.

To find the diagonal entries, observe that both $U_M$ and $UU_M$ are
symmetric, the latter being a scalar multiple of $U_{g_XM}$.
Transposing the product gives $UU_M=U_M U^{\mathsf T}$. Let $w$ be a unit
eigenvector of $U$, with eigenvalue $\lambda$. Conjugating
$U^*w=\bar\lambda w$ gives $U^{\mathsf T}\bar w=\lambda\bar w$, and hence
\[
 U(U_M\bar w)=U_M U^{\mathsf T}\bar w=\lambda U_M\bar w .
\]
Thus $U_M\bar w$ is a unit multiple of $w$, and $U_M^*w$ is a unit multiple of
$\bar w$. The scalar cancels when we take the diagonal entry, giving
\[
 w^*T_M(A)w=(\bar w)^*A^{\mathsf T}\bar w
          =(w^*Aw)^{\mathsf T}=w^*Aw .
\]
Every summand therefore has the same diagonal as $A$ in this eigenbasis.

It remains to count the terms. Changing the sign of the second column of
$M$ gives
\[
 R=\begin{pmatrix}a&-b\\c&-a\end{pmatrix},\qquad
 \operatorname{tr}(RX)=2ar+cs-bt .
\]
As $M$ runs over $\mathcal M$, $R$ runs once over all matrices of trace
zero and determinant $-1$. These two conditions and $\operatorname{tr}(RX)=0$
are invariant under a simultaneous change of basis for $R$ and $X$.
Choose $e,Xe$ independent as in the proof
of Lemma~\ref{lem:groups-bases}. In this basis $X^2=\Delta I$ gives
\[
 X=\begin{pmatrix}0&\Delta\\1&0\end{pmatrix},\qquad
 R=\begin{pmatrix}u&-\Delta v\\v&-u\end{pmatrix},
 \qquad u^2-\Delta v^2=1 ,
\]
where the form of $R$ follows from $\operatorname{tr}(RX)=0$.
The last equation also counts the elements of $C(X)$ in
\eqref{eq:explicit-centralizer}. Thus the sum in \eqref{eq:matrix-local}
has $2c_X$ terms and equals $2c_X E_X(A)$.

\Needspace{9\baselineskip}
\smallskip
\noindent\emph{Counting the classes $X$.}
Now fix $M$ and vary $X$. The associated matrix $R$ has trace zero
and determinant $-1$, so its eigenvalues are $1,-1$.
In an eigenbasis of $R$, the equations $\operatorname{tr}X=\operatorname{tr}(RX)=0$
say that $X$ has the form $\left(\begin{smallmatrix}0&u\\v&0\end{smallmatrix}\right)$.
Of the $p+1$ nonzero pairs $(u,v)$ up to scale, two have $uv=0$.
The others have representatives $(1,v)$ with $v\ne0$. Half of these
have $v$ a square and half a nonsquare. Since $-\det X=uv$, each $M$
occurs equally often in \eqref{eq:matrix-local} for a given type of $X$.
Summing over that type and using \eqref{eq:matrix-total}, we find that
$\sum_{X\text{ of type }a} E_X(A)$ is a constant multiple of
$A+\operatorname{tr}(A)I$, since $c_X$ is constant on each type.
At $A=I$ the sum is $b_a I$ and $A+\operatorname{tr}(A)I=(m+1)I$,
giving the coefficient $b_a/(m+1)$ in \eqref{eq:second-group-sum}.

For the product, multiply \eqref{eq:matrix-local} for $A$ and $B$, and sum
over $X$. Fix $M=\left(\begin{smallmatrix}a&b\\c&a\end{smallmatrix}\right)$
and $N=\left(\begin{smallmatrix}a'&b'\\c'&a'\end{smallmatrix}\right)$.
The term $T_M(A)T_N(B)$ occurs once for each nonzero triple $(r,s,t)$,
up to scale, satisfying
\[
 2ar+cs-bt=0,\qquad 2a'r+c's-b't=0 .
\]
If the equations are independent, they determine one such triple up to scale.
Their coefficients are proportional precisely when $N$ is a scalar multiple
of $M$. As both determinants are one, this means $N=M$ or $N=-M$.
In either case only one equation remains and it gives $p+1$ classes.

We therefore count every ordered pair $M,N$ once, with $p$ further
occurrences for each of $(M,M)$ and $(M,-M)$. On $H$, $U_{-I}=-I$,
so $T_{-M}=T_M$. As transposition reverses the order of the factors,
$T_M(A)T_M(B)=T_M(BA)$. Hence
\[
 4\sum_X c_X^2 E_X(A)E_X(B)
 =\Bigl(\sum_{M\in\mathcal M}T_M(A)\Bigr)
  \Bigl(\sum_{N\in\mathcal M}T_N(B)\Bigr)
  +2p\sum_{M\in\mathcal M}T_M(BA).
\]
Substituting \eqref{eq:matrix-total} gives
\eqref{eq:third-group-sum}, since $\operatorname{tr}(BA)=\operatorname{tr}(AB)$.
\end{proof}

\Needspace{15\baselineskip}

\end{document}